\documentclass{amsart}

\input xy
\xyoption{all}

\begin{document}

\title{On the kernel of the Gassner representation}

\author{Kevin P.~Knudson}\thanks{Partially supported by NSF grant no.~DMS-0242906 and by ORAU}

\address{Department of Mathematics and Statistics, Mississippi State University, Mississippi State, MS 39762}
\email{knudson@math.msstate.edu}

\date{November 1, 2004}


\newtheorem{theorem}{Theorem}[section]
\newtheorem{prop}[theorem]{Proposition}
\newtheorem{lemma}[theorem]{Lemma}
\newtheorem{cor}[theorem]{Corollary}
\newtheorem{conj}[theorem]{Conjecture}
\newtheorem{definition}[theorem]{Definition}
\newtheorem{remark}[theorem]{Remark}

\newcommand{\zz}{{\mathbb Z}}
\newcommand{\zq}{{\mathbb Q}}
\newcommand{\ra}{\rightarrow}
\newcommand{\lra}{\longrightarrow}
\newcommand{\bop}{\bigoplus}
\newcommand{\glnzt}{GL_n(\zz[t_1^{\pm 1},\dots ,t_n^{\pm 1}])}
\newcommand{\lglnz}{M_n(\zz)}
\newcommand{\lslnz}{M^0_n(\zz)}
\newcommand{\gr}{\text{gr}^\bullet}

\begin{abstract}
We study the Gassner representation of the pure braid group $P_n$ by considering its
restriction to a free subgroup $F$.  The kernel of the restriction is shown to lie in the
subgroup $[\Gamma^3F,\Gamma^2F]$, sharpening a result of Lipschutz.
\end{abstract}

\maketitle

\section{Introduction}\label{intro}
Denote by $G_n:P_n\ra\glnzt$ the unreduced Gassner representation of the pure braid group
$P_n$ (a formula is given in
Section \ref{gassner} below).  The faithfulness of $G_n$ for $n\ge 4$ is a long-standing open
question.  In this note, we investigate this by considering the restriction of $G_n$ to a certain
free subgroup $F_{n-1}$ of $P_n$:
$$g_n:F_{n-1}\lra\glnzt.$$
The faithfulness of $G_n$ would follow from that of $g_n$ (see
Proposition \ref{equivalent} below, or \cite{cohen} for a more
general result).

For a group $H$, denote by $\Gamma^\bullet H$ the lower central series of $H$.  The main result
of this paper is the following.

\medskip

\noindent {\bf Theorem \ref{mainthm}.} {\em The kernel of $g_n$ lies in the subgroup
$[\Gamma^3 F_{n-1},\Gamma^2 F_{n-1}]$.}

\medskip

This is proved by passing to the graded quotients associated to
the lower central series of $F_{n-1}$ and the filtration of
$\glnzt$ by powers of the augmentation ideal $J=\ker\{\zz[t_1^{\pm
1}, \dots ,t_n^{\pm 1}]\stackrel{t_i\mapsto 1}{\lra}\zz\}$.  This
allows us to show that the kernel of $g_n$ lies in
$\Gamma^5F_{n-1}$. In \cite{lipschutz}, S.~Lipschutz proved that
the kernel of $g_n$ lies in $[\Gamma^2 F_{n-1},\Gamma^2 F_{n-1}]$
using different techniques (see also \cite{abdul} for another
proof).  These two facts together allow us to prove Theorem
\ref{mainthm}.

  We also show (Theorem \ref{ithterm}) that the intersection of the
kernel of $g_n$ with $\Gamma^sF_{n-1}$ lies in the subgroup
$[\Gamma^{s-2}F_{n-1},\Gamma^2F_{n-1}]\cdot\Gamma^{s+1}F_{n-1}$.

\medskip

\noindent {\em Acknowledgements.}  I thank Fred Cohen for
suggesting this question to me and for many useful conversations.
I am also grateful to an anonymous referee for many helpful
comments.

\section{Preliminaries and Notation}\label{gassner}

\subsection{The Gassner representation} Denote by $A_{rs}$, $1\le r<s\le n$, the generators of $P_n$.
The (unreduced) Gassner representation is the homomorphism
$G_n:P_n\ra\glnzt$ given by the formula:

 \[
G_n(A_{rs}) = \begin{pmatrix}
I_{r-1}&0&0&0&0\\
0        &1-t_{r}+t_{r}t_{s}&0&t_{r}(1-t_{r})&0\\
0        &\vec u      &I_{s-r-1}&\vec v&0\\
0        &1-t_{s}       &0        &t_{r}&0\\
0        &0           &0        &0  &I_{n-s}
\end{pmatrix}
\]
where
\[
\vec u = \begin{pmatrix} (1-t_{r+1})(1-t_{s})&\cdots&
(1-t_{s-1})(1-t_{s}) \end{pmatrix}^\top
\]
and
\[
\vec v = \begin{pmatrix} (1-t_{r+1})(t_{r}-1)&\cdots&
(1-t_{s-1})(t_{r}-1) \end{pmatrix}^\top
\]
and $I_k$ denotes the $k\times k$ identity matrix.  This representation is reducible to
an $(n-1)$-dimensional representation, but the resulting formula is more complicated.

\subsection{The free subgroup}\label{free}
Denote by
$F_{n-1}$ the free subgroup of $P_n$ obtained by deleting the last string; this subgroup has
generators $A_{1n}, A_{2n}, \dots ,A_{n-1,n}$. Moreover, we have a split short exact sequence
$$1\lra F_{n-1} \lra P_n\lra P_{n-1}\lra 1$$
so that $P_n$ is the semidirect product of $P_{n-1}$ and $F_{n-1}$.  Also, the following
diagram commutes (\cite{birman}, p.~138):
$$\xymatrix{
P_n \ar[r]\ar[d]_{G_n} & P_{n-1}\ar[d]^{G_{n-1}} \\
G_n(P_n)\ar[r] & G_{n-1}(P_{n-1}) }$$ where the lower horizontal
map is given by setting $t_n=1$ and deleting the $n$th row and
column.

Denote by $\Gamma^\bullet F_{n-1}$ the lower central series of
$F_{n-1}$ and  for each $i$, consider the free abelian group
$$\Gamma^i F_{n-1}/\Gamma^{i+1}F_{n-1}.$$
  We shall
need an explicit basis of each $\Gamma^i F_{n-1}/\Gamma^{i+1} F_{n-1}$; this is given by the set
of basic commutators of weight $i$.  These are defined as follows.  Denote by $x_j$ the image of
$A_{jn}$ in $F_{n-1}/\Gamma^2 F_{n-1}$.  Then the $x_j$ are the basic commutators of weight one (denote
this by $w(x_j) = 1$) and
having defined the basic commutators of weight less than $i$, the basic commutators of weight $i$ are
the various $[c_u,c_v]$ where

\begin{enumerate}
\item $c_u$ and $c_v$ are basic with $w(c_u)+w(c_v)=i$, and

\item $c_u>c_v$ and if $c_u=[c_a,c_b]$, then $c_v\ge c_b$.
\end{enumerate}
The commutators are ordered as follows.  Those of weight $i$ follow those of weight less than $i$ and are
ordered arbitrarily with respect to each other.  A proof that the basic commutators of weight $i$ form a basis
of $\Gamma^i F_{n-1}/\Gamma^{i+1} F_{n-1}$ may be found in \cite{hall}, p.~175.

Denote by $g_n$ the restriction of $G_n$ to
$F_{n-1}$ and set $X_n=g_n(F_{n-1})$.

\begin{prop}\label{equivalent}
 \begin{eqnarray*}
G_n\,\text{is faithful} & \Leftrightarrow & g_n\,\text{is faithful} \\
                        & \Leftrightarrow & \text{the
                        map}\;\Gamma^iF_{n-1}/\Gamma^{i+1}F_{n-1}\to\Gamma^iX_n/\Gamma^{i+1}X_n
                        \\
                        &   & \text{is injective for each}\;i\ge 1.
\end{eqnarray*}
\end{prop}

\begin{proof} A proof of the first equivalence may be found in, for example, \cite{cohen}.
The second equivalence is an easy exercise about free groups and
is left to the reader.
\end{proof}

\subsection{The congruence subgroup}\label{congruence}
Denote the subgroup of $\glnzt$ consisting of those matrices $A$ with $A\equiv I_n$ modulo
$(t_1-1,t_2-1,\dots ,t_n-1)$ by $K_n$ (i.e., $K_n=GL_n(R,J)$ for $R=\zz[t_1^{\pm 1},\dots ,t_n^{\pm 1}]$ and $J$
the augmentation ideal).  Note that the image of $P_n$ under $G_n$ lies in $K_n$. The group $K_n$ is filtered by powers of $J$:
$$K_n^i = \{A\in K_n:A\equiv I_n\mod J^i\}.$$  This is a central series: $[K_n^i,K_n^j]\subseteq K_n^{i+j}$.

Consider the graded quotients $K_n^i/K_n^{i+1}.$ Note that
$\Gamma^i K_n\subseteq K_n^i$, but it is probably not true that
$K_n^\bullet$ is the lower central series (for $R=\zz[t,t^{-1}]$,
the corresponding group $K_n$ has $K_n^i/\Gamma^iK_n$ a torsion
group for $n\ge 4$ \cite{knudson}).    Consider the induced map
$$\Phi^i: \Gamma^iF_{n-1}/\Gamma^{i+1}F_{n-1}\lra K^i_n/K^{i+1}_n.$$
  Then by Proposition \ref{equivalent}, we have the following:

\medskip

\centerline{$G_n$ is injective if $\Phi^i$ is injective for all
$i\ge 1$.}

\medskip

We show in Section \ref{main} that $\Phi^k$ is injective for $k\le 4$, but that injectivity fails
for $k=5$.

\subsection{Structure of $K^i_n/K^{i+1}_n$}\label{structure}  Given $A\in K_n^i$, we may write
$$A\equiv I_n + \sum_{1\le \ell_1\le\cdots\le \ell_i\le n} (t_{\ell_1}-1)\cdots (t_{\ell_i}-1) A_{\ell_1,\dots ,\ell_i} \mod J^{i+1},$$
where $A_{\ell_1,\dots ,\ell_i}\in \lglnz$.  Define $$\pi_i:K_n^i\lra \bop_{1\le\ell_1\le\cdots\le\ell_i\le n} \lglnz$$ by
$$\pi_i(A) = (A_{\ell_1,\dots ,\ell_i})_{1\le\ell_1\le\dots\le\ell_i\le n}.$$
$\pi_i$ is clearly a homomorphism and $\text{ker}\pi_i=K_n^{i+1}$.

Denote by $e_{\ell m}(a)$ the matrix having $a$ in the $\ell,m$ position and zeroes elsewhere.  Note that
$\pi_1$ is surjective:
$$\pi_1(I_n+e_{\ell m}(t_j-1))=(0,\dots ,0, e_{\ell m}(1),  0,\dots ,0)$$
where $e_{\ell m}(1)$ appears in the summand corresponding to
$(t_j-1)$. Note that this works for $\ell=m$ as $1+(t_j-1)=t_j$ is
a unit in $\zz[t_1^{\pm 1},\dots ,t_n^{\pm 1}]$. For $i>1$, the
image of $\pi_i$ is the sum of copies of $\lslnz$ (matrices of
trace $0$):
$$\pi_i(I_n+e_{\ell m}((t_{j_1}-1)\cdots (t_{j_i}-1)))=
                                    (0,\dots ,0,  e_{\ell m}(1),  0,\dots ,0)$$
occurring in the summand corresponding to the monomial  $(t_{j_1}-1)\cdots (t_{j_i}-1)$
for $\ell\ne m$.  Also, we can hit $e_{\ell,\ell}(1)-e_{\ell+1,\ell+1}(1)$ since
$${\SMALL
U=\left[\begin{array}{cccccccc}
1 & & & & & & &  \\
 & \ddots & & & &  & & \\
 &  & 1 & & & & &  \\
 & & & 1+(t_{j_1}-1)\cdots (t_{j_i}-1) & -(t_{j_1}-1)\cdots (t_{j_i}-1) & & &  \\
 & & & (t_{j_1}-1)\cdots (t_{j_i}-1) & 1-(t_{j_1}-1)\cdots (t_{j_i}-1) & & &  \\
 & & &                               &                                 & 1 & &  \\
 & & &                               &                                 &   & \ddots & \\
 & & &                               &                                 &   &  & 1
 \end{array}\right]\in K_n^i}$$ and $\pi_i(U)=e_{\ell,\ell}(1)-e_{\ell+1,\ell+1}(1)+e_{\ell+1,\ell}(1)-e_{\ell,\ell+1}(1)$.

 \section{The Main Result}\label{main}
  Let us now investigate the map $\Phi^i:\Gamma^iF_{n-1}/\Gamma^{i+1}F_{n-1}\ra  K^i_n/K^{i+1}_n$.
 When we refer to a ``factor" we mean the copy of $\lglnz$ or $\lslnz$ in $K^i_n/K^{i+1}_n$
  corresponding to a certain monomial
 $(t_{j_1}-1)\cdots (t_{j_i}-1)$; we will abbreviate this monomial to $t_{j_1}\cdots t_{j_i}$.

 Now, we have
 $$\begin{array}{rcc}
 \Phi^1(x_r) = & (e_{rr}(1) + e_{nr}(-1), & e_{rn}(-1)+e_{nn}(1)) \\
               &           t_n             &  t_r
                          \end{array}$$
 where the monomial under an entry indicates the factor in which it lies.  For $1\le r\le n-1$, these elements
 are clearly linearly independent in $K_n/K^2_n=\lglnz$, and so $\Phi^1$ is injective.

 On the second level we have, for $r>s$,
 $$\begin{array}{rccc}
 \Phi^2([x_r,x_s]) = & (e_{sr}(-1)+e_{nr}(1), & e_{rs}(1)+e_{ns}(-1), & e_{rn}(-1)+e_{sn}(1)) \\
                     &       t_st_n         & t_rt_n                & t_rt_s
                            \end{array}$$
 and these are clearly linearly independent in $K^2_n/K^3_n$.  So $\Phi^2$ is injective as well.

 To calculate $\Phi^3$, we must order the bases of $F_{n-1}/\Gamma^2F_{n-1}$ and $\Gamma^2F_{n-1}/\Gamma^3F_{n-1}$.
   Use the obvious
 order on the first level:  $x_{n-1}>x_{n-2}>\cdots >x_1$.  On the second level, use
 $$[x_{n-1},x_{n-2}]>\cdots >[x_{n-1},x_1]>[x_{n-2},x_{n-3}]>\cdots >[x_3,x_2]>[x_3,x_1]>[x_2,x_1].$$
 Then a basis of $\Gamma^3F_{n-1}/\Gamma^4F_{n-1}$ is the set
 $$[[x_r,x_s],x_u]\qquad r>s, u\ge s.$$  We have the following formula for $\Phi^3([[x_r,x_s],x_u])$:
 $$\begin{array}{c|l}
 \text{factor} & \text{element} \\ \hline
 t_st_n^2 & e_{sr}(\delta_{us})+e_{sr}(-\delta_{ur})+e_{nr}(-\delta_{us})+e_{nr}(\delta_{ur}) \\
 t_st_ut_n & e_{sn}(\delta_{ru})+e_{ur}(1)+e_{nn}(-\delta_{ru})+e_{nr}(-1) \\
 t_rt_n^2 & e_{rs}(\delta_{us})+e_{rs}(-\delta_{ur})+e_{ns}(\delta_{ur})+e_{ns}(-\delta_{us}) \\
 t_rt_ut_n & e_{rn}(-\delta_{us}) + e_{us}(-1)+e_{nn}(\delta_{us})+e_{ns}(1) \\
 t_st_rt_n & e_{rn}(\delta_{ur})+e_{ru}(1)+e_{nn}(-\delta_{ur})+e_{sn}(-\delta_{us})+e_{su}(-1)+e_{nn}(\delta_{us}) \\
 t_rt_st_u & e_{rn}(-1)+e_{sn}(1)
 \end{array}$$

 Write $c_{rsu}$ for $\Phi^3([[x_r,x_s],x_u])$.

 \begin{prop}\label{level3} $\{c_{rsu}: r>s, u\ge s\}$ is a linearly independent set in $K^3_n/K^4_n$.
 \end{prop}

 \begin{proof}   Suppose
 $$\sum m_{rsu}c_{rsu} =0$$ for some $m_{rsu}\in \zz$.  If $s=u$, the factor
 $t_s^2t_r = t_rt_st_u$ comes into play and does not occur in any other $c_{rsu}$.  So
 $m_{rsu}=0$ in this case.  Similarly if $r=u$, the factor $t_r^2t_s = t_rt_st_u$ comes into
 play and does not occur in any other $c_{rsu}$ and so $m_{rsu}=0$ here as well.

 Thus, we may assume we have $\sum m_{rsu}c_{rsu}=0$ where each of the $c_{rsu}$ has distinct
 $r,s,u$.  For a given fixed $r,s,u$, the factor $t_rt_st_u$ occurs exactly twice---in
 $c_{rsu}$ and $c_{usr}$.  The corresponding elements are $e_{rn}(-1)+e_{sn}(1)$ and
 $e_{un}(-1)+e_{sn}(1)$, respectively.  As $r\ne u$, these are linearly independent in
 $\lglnz$, and so we must have $m_{rsu}=m_{usr}=0$ in this case as well.
 \end{proof}

 Now let's look at $\Phi^4$.  A basis of $\Gamma^4F_{n-1}/\Gamma^5F_{n-1}$ consists of the elements
 $$[[[x_r,x_s],x_u],x_v]\qquad r>s,u\ge s,v\ge u$$
 and
 $$[[x_r,x_s],[x_u,x_v]]\qquad r>s,u>v,r\ge u$$
 (and in addition, if $r=u$, then $s>v$).  Denote the image of an element above under $\Phi^4$ by
 $c_{rsuv}$.  Note that the order $r,s,u,v$ uniquely determines which of the elements we have, as
 no sequence from the first type of basis element can occur as a sequence from the second type.  Under
 $\Phi^4$, the image of $[[[x_r,x_s],x_u],x_v]$ is
 $$\begin{array}{c|l}
 \text{factor} & \text{element} \\ \hline
 t_st_n^3 & e_{sv}(\delta_{us}\delta_{vr})+e_{vr}(-\delta_{vs})+e_{nr}(\delta_{su}\delta_{vs})+e_{sv}(-\delta_{ur}\delta_{vr}) \\
          & {}+e_{vr}(\delta_{ur}\delta_{vs})+e_{nr}(-\delta_{ur}\delta_{vs})+e_{nv}(-\delta_{us}\delta_{vr})+e_{nv}(-\delta_{ur}\delta_{vr}) \\ \hline
 t_st_vt_n^2 & e_{sn}(-\delta_{us}\delta_{rv})+e_{sn}(\delta_{ur}\delta_{rv})+e_{vr}(-\delta_{us})+e_{nn}(\delta_{us}\delta_{rv}) \\
             & {}+e_{nr}(\delta_{us})+e_{vr}(\delta_{ur})+e_{nn}(-\delta_{ur}\delta_{rv})+e_{nr}(-\delta_{ur}) \\ \hline
 t_ut_st_n^2 & e_{vn}(-\delta_{ru}\delta_{vs})+e_{sv}(-\delta_{ru})+e_{nn}(\delta_{ru}\delta_{sv})+e_{uv}(\delta_{rv}) \\
             & {}+e_{vr}(-\delta_{uv})+e_{nr}(\delta_{uv})+e_{nv}(\delta_{ru})+e_{nv}(-\delta_{rv}) \\ \hline
 t_ut_st_vt_n & e_{sn}(\delta_{ru})+e_{un}(-\delta_{rv})+e_{vn}(-\delta_{ru}) \\
              & {}+ e_{vr}(-1) + e_{nn}(\delta_{rv})+e_{nr}(1) \\ \hline
 t_rt_n^3 & e_{rv}(\delta_{us}\delta_{vs})+e_{vs}(-\delta_{us}\delta_{vr})+e_{ns}(\delta_{us}\delta_{vr})+e_{rv}(-\delta_{ur}\delta_{vs}) \\
          & {}+e_{vs}(\delta_{ur}\delta_{vr})+e_{ns}(\delta_{ur}\delta_{vr})+e_{nv}(\delta_{ur}\delta_{sv})+e_{nv}(-\delta_{us}\delta_{vs}) \\ \hline
 t_rt_vt_n^2 & e_{rn}(-\delta_{us}\delta_{vs})+e_{rn}(\delta_{ur}\delta_{vs})+e_{vs}(\delta_{ur})+e_{nn}(-\delta_{ur}\delta_{vs}) \\
             & {}+e_{ns}(-\delta_{ur})+e_{vs}(-\delta_{us})+e_{nn}(\delta_{us}\delta_{vs})+e_{ns}(\delta_{us}) \\ \hline
 t_rt_ut_n^2 & e_{vn}(\delta_{us}\delta_{vr})+e_{rv}(\delta_{us})+e_{nn}(-\delta_{us}\delta_{vr})+e_{uv}(-\delta_{vs}) \\
             & {}+e_{vs}(\delta_{uv})+e_{ns}(-\delta_{uv})+e_{nv}(-\delta_{us})+e_{nv}(\delta_{vs}) \\ \hline
 t_rt_ut_vt_n & e_{rn}(-\delta_{us})+e_{un}(\delta_{vs})+e_{vn}(\delta_{us}) \\
              & {}+e_{vs}(1)+e_{nn}(-\delta_{vs})+e_{ns}(-1) \\ \hline
 t_rt_st_n^2 & e_{vn}(-\delta_{ur}\delta_{vr})+e_{rv}(-\delta_{ur})+e_{nn}(\delta_{ur}\delta_{vr})+e_{rv}(\delta_{uv})+e_{vu}(-\delta_{vr}) \\
             & {}+e_{nu}(\delta_{vr})+e_{nv}(\delta_{ur})+e_{vn}(\delta_{us}\delta_{vs})+e_{sv}(\delta_{us})+e_{nn}(-\delta_{us}\delta_{vs}) \\
             & {}+e_{sv}(-1)+e_{vu}(\delta_{vs})+e_{nu}(-\delta_{vs})+e_{nu}(-\delta_{vs})+e_{nv}(-\delta_{vs}) \\ \hline
 t_rt_st_vt_n & e_{rn}(\delta_{ur})+e_{rn}(-\delta_{uv})+e_{vn}(-\delta_{ur}) \\
              & {}+e_{sn}(-\delta_{us})+e_{sn}(\delta_{uv})+e_{vn}(\delta_{us}) \\ \hline
 t_rt_st_ut_n & e_{vn}(\delta_{rv})+e_{vn}(-\delta_{sv})+e_{rv}(1) \\
              & {}+e_{nn}(-\delta_{rv})+e_{sv}(-1)+e_{nn}(\delta_{sv}) \\ \hline
 t_rt_st_ut_v & e_{rn}(-1)+e_{sn}(1)
 \end{array}$$
and the image of $[[x_r,x_s],[x_u,x_v]]$ is
$$\begin{array}{c|l}
\text{factor} & \text{element} \\ \hline
t_st_vt_n^2 & e_{su}(\delta_{rv})+e_{vr}(-\delta_{us})+e_{nr}(\delta_{us})+e_{nu}(-\delta_{rv}) \\ \hline
t_st_ut_n^2 & e_{sv}(-\delta_{ru})+e_{ur}(\delta_{vs})+e_{nr}(-\delta_{vs})+e_{nv}(\delta_{ru}) \\ \hline
t_st_ut_vt_n & e_{sn}(\delta_{ru})+e_{sn}(\delta_{rv})+e_{ur}(1) \\
             & {}+e_{nn}(-\delta_{ru})+e_{vr}(-1)+e_{nn}(\delta_{rv}) \\ \hline
t_rt_vt_n^2 & e_{ru}(-\delta_{sv})+e_{vs}(\delta_{ru})+e_{ns}(-\delta_{ru})+e_{nu}(\delta_{sv}) \\ \hline
t_rt_ut_n^2 & e_{rv}(\delta_{us})+e_{us}(-\delta_{rv})+e_{ns}(\delta_{rv})+e_{nv}(-\delta_{su}) \\ \hline
t_rt_ut_vt_n & e_{rn}(-\delta_{us})+e_{rn}(\delta_{vs})+e_{us}(-1) \\
             & {}+ e_{nn}(\delta_{us})+e_{vs}(1)+e_{nn}(-\delta_{vs}) \\ \hline
t_rt_st_vt_n & e_{vn}(-\delta_{ru})+e_{ru}(-1)+e_{nn}(\delta_{ru}) \\
             & {}+e_{vn}(\delta_{su})+e_{su}(1)+e_{nn}(-\delta_{su}) \\ \hline
t_rt_st_ut_n & e_{un}(\delta_{rv})+e_{rv}(1)+e_{nn}(-\delta_{rv}) \\
             & {}+e_{un}(-\delta_{sv})+e_{sv}(-1)+e_{nn}(\delta_{sv}) \\ \hline
t_rt_st_ut_v & 0
\end{array}$$

\begin{prop}\label{level4} The elements $\{c_{rsuv}\}$ are linearly independent in $K^4_n/K^5_n$.
\end{prop}

\begin{proof} Suppose $\sum m_{rsuv}c_{rsuv}=0$ for some $m_{rsuv}\in \zz$.  Consider the case where
only $2$ of the $r,s,u,v$ are distinct (e.g. $c_{2111}$, $c_{3133}$, etc.).  Then $c_{rsuv}$ contributes to
the factor $t_rt_st_ut_v=t_i^kt_j^\ell$ where $i$ and $j$ are the distinct indices and $1\le k,\ell\le 3$,
$k+\ell=4$.  Note that $c_{rsuv}$ is the only contributor to this factor as the choice of $r$ determines the sequence---the
only possibilities are $rsrr$, $rsss$ or $rssr$ (note that no $[[x_r,x_s],[x_u,x_v]]$ occur as $r\ne s$, $u\ne v$ implies that
$r=u$ and $s=v$ and so the element is $0$).  To this factor, $c_{rsuv}$ contributes $e_{rn}(-1)+e_{sn}(1)$.  It follows
that $m_{rsuv}=0$ for these elements.

Now suppose that $c_{rsuv}=\Phi([[[x_r,x_s],x_u],x_v])$ with $r,s,u,v$ distinct (we shall deal with the double commutators
with 4 distinct indices below).  Then $c_{rsuv}$ contributes the element $e_{rn}(-1)+e_{sn}(1)$ to the factor
$t_rt_st_ut_v$.  For fixed $r,s,u,v$, we must have $r>s$ and $s<u<v$ if they are all distinct.  So the only contributors to
this factor are (1) $c_{rsuv}$, (2) $c_{usrv}$ or $c_{usvr}$ (the latter if $v<r$), (3) $c_{vsru}$ or $c_{vsur}$ (the latter if
$u<r$); and they contribute $e_{rn}(-1)+e_{sn}(1)$, $e_{un}(-1)+e_{sn}(1)$, $e_{vn}(-1)+e_{sn}(1)$ respectively.  Since these
elements are linearly independent in $\lglnz$, we must have $m_{rsuv}=m_{usrv}=m_{vsru}=0$.

Next, suppose that $r,s,u,v$ consist of 3 distinct indices.  We have $r>s$, $u\ge s$, and $v\ge u$.  There are three cases to consider.
Recall that in any case, $r>s$.

\medskip

\noindent {\bf Case 1.}    $s=u$.  We may assume (without loss of generality) that $r>v$ and $v>u$.  Then we have
three possible elements to consider:  $c_{rssv}$, $c_{vssr}$, and $c_{rsvs}$ (the latter corresponds to $[[x_r,x_s],[x_v,x_s]]$).
Here, the factor $t_s^2t_vt_r$ receives contributions only from $c_{rssv}$ and $c_{vssr}$, the elements being
$e_{rn}(-1)+e_{sn}(1)$ and $e_{vn}(-1)+e_{sn}(1)$, respectively.  These are linearly independent in $\lglnz$ and so
$m_{rssv}=m_{vssr}=0$.  Then, in the factor $t_s^2t_rt_n$, the only contributors are $c_{rssv}$, $c_{vssr}$ and
$c_{rsvs}$---the latter contributing $e_{rv}(1)+e_{sv}(-1)$, while the former two contribute $e_{rv}(1)+e_{sv}(-1)$ and
$e_{vr}(1)+e_{sr}(-1)$.  As we've already shown that $m_{rssv}=m_{vssr}=0$, and since the only other $c_{rsuv}$ that contribute
here only have two distinct indices, we must have $m_{rsvs}=0$ as well.

\medskip

\noindent {\bf Case 2.}  $u=v$.  Then $u>s$.  Suppose $r>u$ (the case $r<u$ is Case 3 below).  Then the three elements to
consider are $c_{rsuu}$, $c_{usur}$ and $c_{ruus}$.  The first two elements contribute to the factor $t_st_u^2t_r$ the elements
$e_{rn}(-1)+e_{sn}(1)$ and $e_{un}(-1)+e_{sn}(1)$ respectively and no other $c_{ijkl}$ contributes to this factor.  So $m_{rsuu}=m_{usur}=0$.
Then consider the factor $t_u^2t_rt_n$.  Here, $c_{ruus}$ contributes $e_{ru}(1)+e_{su}(-1)$ and the only other contributors
have $m_{ijkl}=0$ already.  Thus, $m_{ruus}=0$, as well.

\medskip

\noindent {\bf Case 3.}  $u=v$, $u>r$.  This is similar to Case 2.

\medskip

Finally, consider the $c_{rsuv}=\Phi([[x_r,x_s],[x_u,x_v]])$ with $r,s,u,v$ distinct.  Then
$r>s$, $u>v$, $r>u$.  As the indices are distinct the only factors contributed to are
$t_st_ut_vt_n$, $t_rt_ut_vt_n$, $t_rt_st_vt_n$, and $t_rt_st_ut_n$ (see the formulas above).  For
fixed $r,s,u,v$, the only possible $c_{ijkl}$ are given in the following table, along with the elements
contributed to each factor.
{\small $$\begin{array}{r|c|c|c|c}
          & t_st_ut_vt_n & t_rt_ut_vt_n & t_rt_st_vt_n & t_rt_st_ut_n \\ \hline
c_{rsuv}  & e_{ur}(1)+e_{vr}(-1) & e_{us}(-1)+e_{vs}(1) & e_{ru}(-1)+e_{su}(1) & e_{rv}(1)+e_{sv}(-1) \\ \hline
c_{rusv}  & e_{sr}(1)+e_{vr}(-1) & e_{rs}(-1)+e_{us}(1) & e_{su}(-1)+e_{vu}(1) & e_{rv}(1)+e_{uv}(-1) \\ \hline
c_{rvsu}  & e_{sr}(1)+e_{ur}(-1) & e_{rs}(-1)+e_{vs}(1) & e_{ru}(1)+e_{vu}(-1) & e_{sv}(-1)+e_{uv}(1)
\end{array}$$}
In each factor, we obtain linearly dependent elements, but we must remember that we're scaling the element coming from
$c_{rsuv}$ by $m_{rsuv}$.  Looking at the factor $t_st_ut_vt_n$, we find that $m_{rsuv}=m_{rvsu}$ and $m_{rusv}=-m_{rvsu}$.
But then looking at the factor $t_rt_ut_vt_n$, we find $m_{rsuv}=-m_{rvsu}$ and $m_{rusv}=-m_{rvsu}$.  Thus,
$m_{rsuv}=m_{rusv}=m_{rvsu}=0$.

This completes the proof.
\end{proof}

\begin{cor}\label{kernelgamma5} The kernel of $g_n$ is contained in $\Gamma^5F_{n-1}$.\hfill $\qed$
\end{cor}

It is possible to sharpen Corollary \ref{kernelgamma5} to obtain the main result.

\begin{theorem}\label{mainthm} The kernel of $g_n$ is contained in $[\Gamma^3F_{n-1},\Gamma^2F_{n-1}]$.
\end{theorem}

\begin{proof} For simplicity, denote the group $\Gamma^iF_{n-1}$ by $\Gamma^i$.
By Corollary \ref{kernelgamma5} and by \cite{lipschutz}, we have
$\text{ker}(g_n)\subseteq \Gamma^5\cap [\Gamma^2,\Gamma^2]$.  We
claim that the latter group equals $[\Gamma^3,\Gamma^2]$.  The
main theorem in \cite{hurley1} implies that
$$\Gamma^5\cap [\Gamma^2,\Gamma^2] =
I_{\Gamma^2}([\Gamma^3,\Gamma^2]),$$ where $I_R(S)$ is the
isolator of $S$ in $R$. (Recall that the isolator of $S$ in $R$ is
the set $I_R(S)=\{x\in R: x^n\in S\,\text{for some}\,n\}$.) To see
that this latter group is simply $[\Gamma^3,\Gamma^2]$, it
suffices to show that the quotient group
$\Gamma^2/[\Gamma^3,\Gamma^2]$ is torsion-free. Consider the short
exact sequence
$$1\lra \frac{[\Gamma^2,\Gamma^2]}{[\Gamma^3,\Gamma^2]} \lra
\frac{\Gamma^2}{[\Gamma^3,\Gamma^2]}\lra\frac{\Gamma^2}{[\Gamma^2,\Gamma^2]}\lra
1.$$ Since $\Gamma^2$ is a free group, the group
$\Gamma^2/[\Gamma^2,\Gamma^2]$ is free abelian. But, by Theorem 6
of \cite{hurley2}, the group
$[\Gamma^2,\Gamma^2]/[\Gamma^3,\Gamma^2]$ is also free abelian. It
follows that $\Gamma^2/[\Gamma^3,\Gamma^2]$ is torsion-free.  This
completes the proof.
\end{proof}

The methods used above allow us to prove the following result.

\begin{theorem}\label{ithterm} For $s\ge 5$, $\text{\em
ker}(g_n)\cap \Gamma^s \subseteq [\Gamma^{s-2},\Gamma^2]\cdot
\Gamma^{s+1}$.
\end{theorem}

\begin{proof}
Note that any basic commutator is given by a unique list of
integers corresponding to the $x_j$ that occur in the commutator.
For example, $[[x_3,x_2],[x_3,x_1]]$ yields the list $3,2,3,1$. We
therefore may denote a basic commutator of weight $s$ by
$x_{\ell_1\ell_2\cdots\ell_s}$ without confusion.  Denote the
element $\Phi^s(x_{\ell_1\ell_2\cdots\ell_s})$ by
$c_{\ell_1\ell_2\cdots\ell_s}$. Note that all the basic
commutators in $\Gamma^s/\Gamma^{s+1}$ lie in
$[\Gamma^{s-2},\Gamma^2]$, except for the various $[c_u,x_j]$.
Moreover, any element of the latter form must be an $s$-fold
commutator:
$$x_{\ell_1\ell_2\cdots\ell_s}=[\cdots [[[x_{\ell_1},x_{\ell_2}],x_{\ell_3}],\cdots ],x_{\ell_s}].$$
 To prove the theorem, it suffices to show that if we
have a dependency relation
$$\sum m_{\ell_1\ell_2\cdots\ell_s} c_{\ell_1\ell_2\cdots\ell_s} = 0$$
where the $m_{\ell_1\cdots\ell_s}\in\zz$, then we have
$m_{\ell_1\cdots\ell_s}=0$ whenever $x_{\ell_1\cdots\ell_s}$ is an
$s$-fold commutator.  This will show that the $[c_u,x_j]$ inject
into $K_n^s/K_n^{s+1}$ and hence that the intersection of the
kernel of $g_n$ with $\Gamma^s$ lies in
$[\Gamma^{s-2},\Gamma^2]\cdot \Gamma^{s+1}$.

Observe that in the case of an $s$-fold commutator, the element
$c_{\ell_1\cdots\ell_s}$ contributes the element
$e_{\ell_1,n}(-1)+e_{\ell_2,n}(1)$ to the factor
$t_{\ell_1}t_{\ell_2}\cdots t_{\ell_s}$ (this is easily proved by
induction using the formulas given above for the $c_j$ and
$c_{rs}$).  Suppose the $\ell_1,\ell_2,\dots ,\ell_s$ consist of
$i$ distinct indices, say $r_1=\ell_1$, $r_2=\ell_2$, and
$r_3,\dots ,r_i$.  We have $r_1>r_2$ and $r_2<r_3<\cdots <r_i$. We
have several contributors to the factor
$t_{\ell_1}t_{\ell_2}\cdots t_{\ell_s} =
t_{r_1}^{a_1}t_{r_2}^{a_2}\cdots t_{r_i}^{a_i}$ (here $a_k$ is the
number of times $r_k$ occurs).  Let us abbreviate notation and
write $c_{r_1r_2\cdots r_i}$ for $c_{\ell_1\cdots\ell_s}$. Certain
permutations of the $\ell_j$ yield $s$-fold basic commutators;
each of these contributes to the factor
$t_{r_1}^{a_1}t_{r_2}^{a_2}\cdots t_{r_i}^{a_i}$ under
consideration. We must show that the resulting contributions are
linearly independent.  Note that the only contributors to this
factor are $s$-fold commutators---if $c=[c_u,c_v]$ where
$w(c_v)\ge 2$, then every factor to which $c$ contributes contains
a power of $t_n$ (see the formulas above).  Thus, we may detect
any dependency relation among the $c_{\ell_1,\dots ,\ell_s}$ by
considering only the factor $t_{\ell_1}\cdots t_{\ell_s}$.

 Now, for some $q$ with $3\le q\le i$ we
must have $r_q<r_1< r_{q+1}$.  Then we get the following contributions to the
factor $t_{r_1}^{a_1}t_{r_2}^{a_2}\cdots t_{r_i}^{a_i}$:
$$\begin{array}{c|c}
\text{sequence}  & \text{element}  \\ \hline
{r_1,r_2,\dots ,r_i} & e_{r_1,n}(-1) + e_{r_2,n}(1)  \\
{r_q,r_2,\dots ,r_{q-1},r_1,r_{q+1},\dots ,r_i}  & e_{r_q,n}(-1)+e_{r_2,n}(1)  \\
{r_{q+1},r_2,\dots ,r_q,r_1,r_{q+2},\cdots ,r_i} & e_{r_{q+1},n}(-1) + e_{r_2,n}(1)  \\
{r_{q+2},r_2,\dots ,r_{q},r_1,r_{q+1},r_{q+3},\dots ,r_i} & e_{r_{q+2},n}(-1) + e_{r_2,n}(1) \\
\vdots  &  \vdots  \\
{r_i,r_2,\dots ,r_q,r_1,r_{q+1},\dots ,r_{i-1}} & e_{r_i,n}(-1) + e_{r_2,n}(1)
\end{array}$$
Since $r_1,r_2,\dots ,r_i$ are distinct, the elements in the second column are linearly independent in
$\lglnz$ (as $i<n$) and so each of the corresponding coefficients satisfies $m_{\ell_1\cdots\ell_s}=0$.
This completes the proof.
\end{proof}

\section{Breakdown}\label{counter}
The method used in Section \ref{main} breaks down at the fifth
level, however.  Indeed, if $n=4$, the kernel of $\Phi^5$ is
rather large.  For example, we have $c_{21131}=c_{31121}$.  This
has the interpretation that the degree $5$ part of the polynomials
in the Gassner matrices of
$$[[[A_{24},A_{14}],A_{14}],[A_{34},A_{14}]]$$ and
$$[[[A_{34},A_{14}],A_{14}],[A_{24},A_{14}]]$$ are the same.  The matrices are not the same,
however, and a computer search
by the author based on the relations in the kernel of $\Phi^5$ has not turned up any elements in the
 kernel of $G_4$.

Note, however, that the failure of the method does not imply that
$g_n$ is not injective.
 Really, one needs
to consider the quotients $\Gamma^i g_n(F_{n-1})/\Gamma^{i+1}g_n(F_{n-1})$ rather
than the classes of the various elements
in $\Gamma^i g_n(F_{n-1})$ modulo the subgroup $K_n^{i+1}$.  This seems to be rather
intractable, however, given the ranks
of the various $\Gamma^iF_{n-1}/\Gamma^{i+1}F_{n-1}$ (for example, $\Gamma^5F_3/\Gamma^6F_3$ has rank $116$).

\end{document}